\documentclass[12pt]{amsart}
\usepackage{amssymb}
\usepackage[dvips]{graphics}
\ifx\mathrm\undefined\let\mathrm\rm\fi
\ifx\mathbf\undefined\let\mathbf\bf\fi
\ifx\mathfrak\undefined\let\mathfrak\frak\fi
\ifx\mathcal\undefined\let\mathcal\cal\fi
\ifx\mathbb\undefined\let\mathbb\Bbb\fi
\ifx\emph\undefined\let\emph\it\fi
\newcommand{\g}{{{\mathfrak g}\,}}
\newcommand{\half}{{\frac12}}

\newcommand{\h}{{{\mathfrak h\,}}}

\newcommand{\Z}{{\mathbb Z}}

\newcommand{\R}{{\mathbb R}}
\newcommand{\C}{{\mathbb C}}

\newcommand{\Ref}[1]{{(\ref{#1})}}
\newcommand{\be}{\begin{displaymath}}
\newcommand{\ee}{\end{displaymath}}
\newcommand{\bea}{\begin{eqnarray*}}
\newcommand{\eea}{\end{eqnarray*}}

\newcommand{\Mu}{{\mathrm{M}}}
\newcommand{\semidirect}{\tilde\times}
\newcommand{\LLambda}{{\vec{\Lambda}}}
\newcommand{\zz}{{\vec{z}}}
\newcommand{\dontprint}[1]{\relax}
\newenvironment{prf}{\noindent{\it Proof\/}:}{$\;\square$
\par\medskip}

\newtheorem%
{thm}{Theorem}[section]
\newtheorem%
{proposition}[thm]{Proposition}
\newtheorem%
{lemma}[thm]{Lemma}
\newtheorem%
{lemmadef}[thm]{Lemma-Definition}
\newtheorem%
{corollary}[thm]{Corollary}
\newtheorem%
{conjecture}[thm]{Conjecture}

\title[{}]{The q-deformed Knizhnik--Zamolodchikov--Bernard
heat equation}
\author[{}]
{Giovanni Felder${}^{*}$ 
\and Alexander Varchenko${}^{**,1}$}
\thanks{${}^1$Supported in part by NSF grant  DMS-9801582}
\begin{document}
\maketitle

\medskip
\centerline{\it ${}^*$Departement Mathematik, ETH-Zentrum,}
\centerline{\it 8092
Z\"urich, Switzerland}
\medskip
\centerline{\it ${}^{**}$Department of Mathematics,
University of North Carolina at Chapel Hill,}
\centerline{\it Chapel Hill, NC 27599-3250, USA}
\medskip
\begin{abstract}
  We introduce a $q$-deformation of the genus one $sl_2$
  Knizhnik--Zamolodchikov--Bernard heat equation. We show that this
  equation for the dependence on the moduli of elliptic curves is
  compatible with the qKZB equations, which give the dependence on the
  marked points.
\end{abstract}
\section{Introduction}
The Knizhnik--Zamolodchikov--Bernard equations are a system of
differential equations
arising in conformal field theory on Riemann surfaces. 
For each $g,n\in \Z_{\geq 0}$, a simple complex Lie algebra $\g$,
$n$ highest weight $\g$-modules $V_i$, and a complex parameter
$\kappa$, we have such a system
of equations. In the case
case of genus $g=1$, they have the form
\begin{equation}\label{eq-KZB}
\kappa\partial_{z_j}v=
-\sum_\nu h_\nu^{(j)}\partial_{\lambda_\nu}v +\sum_{l:l\neq j}
r(z_j-z_l,\lambda)^{(j,l)}v.
\end{equation}
The unknown function $v$ takes values in the zero weight space
$V[0]=\cap_{x\in\h}\mathrm{Ker}(x)$ 
of the tensor product $V=V_1\otimes\cdots\otimes V_n$ with
respect to a  Cartan subalgebra $\h$ of $\g$.
It  depends on variables $z_1,\dots,z_n\in\C$ and $\lambda=\sum\lambda_\nu
h_\nu\in\h$, where $(h_\nu)$ is an orthonormal basis of $\h$, with
respect to a fixed invariant bilinear form.
The notation $x^{(j)}$ for $x\in\mathrm{End}(V_j)$
or $x\in \g$ means $1\otimes\cdots\otimes x
\otimes\cdots\otimes 1$. Similarly $x^{(i,j)}$
denotes the action on the $i$th  and $j$th
factor of $x\in\mathrm{End}(V_i\otimes V_j)$.

 The ``$r$-matrix'' $r\in\g\otimes\g$ obeys
\[
r(z,\lambda)+
r(-z,\lambda)^{(2,1)}=0,\qquad [r(z,\lambda),h\otimes 1+1\otimes h]=0,
\forall h\in\h,
\]
with $(\sum_ix_i\otimes y_i)^{(21)}=\sum_iy_i\otimes x_i$,
and is a solution of the classical dynamical 
Yang--Baxter equation \cite{FW}
($r^{(1,2)}=r(z_1-z_2,\lambda)\otimes 1\in U\g^{\otimes 3}$ etc.)
\begin{eqnarray*}
\sum_\nu
\partial_{\lambda_\nu}r^{(1,2)}h_\nu^{(3)}+
\sum_\nu\partial_{\lambda_\nu}r^{(2,3)}h_\nu^{(1)}+
\sum_\nu\partial_{\lambda_\nu}r^{(3,1)}h_\nu^{(2)} & & \\
-[r^{(1,2)},r^{(1,3)}]
-[r^{(1,2)},r^{(2,3)}]
-[r^{(1,3)},r^{(2,3)}] &=& 0.
\end{eqnarray*}
As a consequence, the KZB equations \Ref{eq-KZB} are compatible, meaning
that if the equations are written as $\nabla_jv=0$, then 
the differential operators  $\nabla_j$ commute with each other.
The solutions of the classical dynamical Yang--Baxter equation
arising in conformal field theory are parametrized by the modulus
$\tau$ in the upper half plane and can be expressed in terms of
theta functions, see \cite{FW,FV1}.

A difference version of this story was proposed in \cite{F}:
Suppose that for an Abelian complex 
Lie algebra $\h$ we have $\h$-modules $V_i$, $i=1,\dots, n$
with a weight decomposition $V_i=\oplus_{\mu\in\h^*}V_i[\mu]$ into
finite dimensional weight spaces $V_i[\mu]$. Then we say that 
meromorphic functions 
$R_{ij}(z,\lambda)$ of $z\in\C$ and $\lambda\in\h^*$
with values in $\mathrm{End}_{\h}(V_i\otimes V_j)$, $(1\leq i\neq j\leq n)$
form a {\em system of dynamical $R$-matrices} if they obey the (quantum)
dynamical Yang--Baxter
equation
\begin{eqnarray*}
\lefteqn{R_{ij}(z_1-z_2,\lambda-2\eta h^{(3)})^{(12)}
R_{ik}(z_1-z_3,\lambda)^{(13)}
R_{jk}(z_2-z_3,\lambda-2\eta h^{(1)})^{(23)}}
\\
 & &
=
R_{jk}(z_2-z_3,\lambda)^{(23)}
R_{ik}(z_1-z_3,\lambda-2\eta h^{(2)})^{(13)}
R_{ij}(z_1-z_2,\lambda)^{(12)},
\end{eqnarray*}
in $\mathrm{End}(V_i\otimes V_j\otimes V_k)$ for all 
$i<j<k$ and are ``unitary'':
\[
R_{ij}(z,\lambda)R_{ji}(-z,\lambda)^{(21)}=\mathrm{Id}_{V_i\otimes V_j}.
\]
We adopt a standard notation: for instance,
$R(z,\lambda-2\eta h^{(3)})^{(12)}$ acts on 
a tensor
$v_1\otimes v_2\otimes v_3$ as 
$R(z,\lambda-2\eta\mu_3)\otimes{\mathrm{Id}}$
if $v_3$ has weight $\mu_3$.

The deformation parameter (``Planck's constant'') is here $\eta$.
If we have a family of dynamical $R$-matrices depending on $\eta$
such that $R_{ij}=\mathrm{Id}_{V_i\otimes V_j}+2\eta r_{ij}+O(\eta^2)$ as $\eta\to 0$, 
we recover the classical dynamical Yang--Baxter equation and
the unitarity condition for $r_{ij}$, i.~e.\ the properties
that $r$ obeys, viewed as an element of $\mathrm{End}(V_i\otimes V_j)$.

If we have a system of 
dynamical $R$-matrices $R_{ij}$ we can then construct a
compatible system of difference equations, the qKZB equations for a
function $v(z_1,\dots,z_n,\lambda)\in V[0]$. They are a dynamical
version of the I. Frenkel--Reshetikhin qKZ equations \cite{FR}, and
their semiclassical limit are the KZB equations. Their construction
is reviewed in \ref{suse-11} below.

The main examples of solutions of the classical and of the
quantum dynamical Yang--Baxter equations are associated with
elliptic curves.
In the quantum case, they can be viewed as intertwining operators
between tensor products of representations of elliptic
quantum groups taken in different orders \cite{FV2}. 
In the rank one case
(one-dimensional $\h$) explicit expressions for $R$ matrices
$R_{\Lambda,\Mu}$ depending on two ``highest weights''
$\Lambda,\Mu\in\C$ are known. 
They are associated to pairs of evaluation Verma modules \cite{FV2}
for the elliptic quantum group $E_{\tau,\eta}(sl_2)$ and were
computed using the functional realization of these modules
\cite{FTV1}. If $n$ highest weights $\Lambda_1,\dots,\Lambda_n\in\C$
are given, then $R_{ij}=R_{\Lambda_i,\Lambda_j}$ form a system
of dynamical $R$ matrices as described above.

Hypergeometric solutions of the corresponding qKZB equations were 
introduced and studied in \cite{FTV1}, \cite{FTV2}, \cite{FV4}.
See also \cite{T} where similar equations are studied and
solved.
Special cases of these equations appear in the statistical
mechanics of RSOS models. Form factors and correlation
functions in the infinite volume limit are conjectured
to obey qKZB equations. In these cases explicit formulae
were proposed by Lukyanov and Pugai \cite{LP}.

The subject of this paper is a deformation of the 
KZB heat equation: in
the classical case, additionally to the KZB equations above,
that are associated to
the variation of the marked points on the elliptic curve, one
also has an equation associated to the variation of
the modulus $\tau$ of the elliptic curve.
The function $v$ also depends on  $\tau$ 
 and one has an additional equation, compatible
with the KZB equations, the
KZB {\em heat equation} 
\[
4\pi i\kappa\partial_\tau v=\triangle_\lambda v
+\frac12\sum_{i,j}s(z,\lambda,\tau)^{(ij)}v.
\]
for some $s\in\g\otimes\g$. Here $\triangle_\lambda$ denotes
the Laplacian of $\h$ corresponding to the invariant bilinear form.
For example, if $n=1$ then
this equation reduces to 
\[
4\pi i\kappa\partial_\tau v=
(\triangle_\lambda - \sum_{\alpha\in\Delta}\wp(\alpha(\lambda),\tau)
e_\alpha e_{-\alpha})v,
\]
where $e_\alpha$ are properly normalized root vectors
and $\wp$ is the Weierstrass function.

In this paper we propose a discrete version of the KZB heat equation
in the rank one case. The heat operator is an integral operator,
whose kernel is given in terms of hypergeometric integral solutions
of the qKZB equations of \cite{FTV2}.

In Sect.\ \ref{se-1} we review the qKZB equations and their 
hypergeometric solutions. Then we introduce the 
elliptic Shapovalov form, which is an ingredient
in the integral
operator, and the qKZB heat equation in Sect.\ \ref{se-2}.
We prove that it is compatible with the qKZB equations, discuss
its properties
and show in Sect.\ \ref{suse-23}, in an illustrative
example, that its semiclassical limit coincides with
the KZB heat equation.
Finally, in Sect.~\ref{se-4} 
we study the special case where the step of the difference
equation is a negative integer multiple of the deformation parameter. In this
case, the semiclassical limit gives the KZB equations with positive integer
$\kappa$, the situation arising in conformal field theory. In this case
the the KZB equations are defined on sections of the finite dimensional
vector bundle of
conformal blocks, of which we describe a difference analogue in simple
cases.

It is very likely that the hypergeometric solutions of the
qKZB equations are also solutions of the qKZB heat equation.
However, we were able to prove this only in the case where
the sum of the highest weights is two. In this case
the hypergeometric integrals are one-dimensional.

In a sequel to this paper, we show that integral operators
of the kind introduced in this paper can also be used to
describe the transformation properties of hypergeometric solutions
under the modular group. In fact it turns out that the hypergeometric
solutions, at least if $\sum\Lambda_i=2$ obey remarkable identities
under transformations of the modulus $\tau$ and the step
$p$  by $\mathrm{SL}(3,\Z)$ acting on $\C P^2$ with affine coordinates
$\tau,p$. These identities give both the solutions and the monodromy
of the solutions.
The whole picture results in an non-Abelian version of the 
properties of the elliptic gamma functions, which is a generalized
Jacobi modular form for $SL(3,\Z)\times\Z^3$ in the sense of \cite{FV5}.

\medskip

\noindent{\bf Acknowledgment.} We thank R. Ferretti for
explanations on Gauss sums.

\section{Hypergeometric solutions of the qKZB equations}\label{se-1}

\subsection{The qKZB equations}\label{suse-11}
Fix $\LLambda=(\Lambda_1,\dots,\Lambda_n)\in\R^n$ such
that $m=\sum_{i=1}^n\Lambda_i/2$ is a nonnegative integer, and a
complex number $\eta$. Unless stated otherwise, we will assume
that these parameters are generic.

Let $\tau$ and $p$ be generic complex numbers in the upper half plane.

Let $V_{\Lambda_j}$ be the vector space with basis $e_0, e_1,\dots$
equipped with the action of an operator $h$ given by 
$he_k=(\Lambda_j-2k)e_k$. We view $V_{\Lambda_i}$ as a
representation of the Abelian Lie  algebra $\h=\C h$. 
 
To these data is associated a system of dynamical $R$-matrices
and thus a system of qKZB difference equations.
The $R$-matrices $R_{\Lambda_j,\Lambda_k}(z,\lambda,\tau)$  \cite{FV2}
of $E_{\tau,\eta}(sl_2)$ are
endomorphisms of $V_{\Lambda_j}\otimes V_{\Lambda_k}$.
Let $V_\LLambda
=V_{\Lambda_1}\otimes\cdots\otimes V_{\Lambda_n}$. The 
qKZB equations are
 equations for
a meromorphic function $v(\zz,\lambda)$ of $\zz\in\C^n$ and
$\lambda\in\C$ taking its values in the zero
weight subspace $V_\LLambda[0]=\mathrm{Ker}(\sum_{i=1}^nh^{(i)})$
of $V_\LLambda$ (this subspace is nontrivial since $\sum\Lambda_i/2$
is assumed to be a nonnegative  integer). 
It will be more convenient to view $v$ as
a function $v(\zz)$ taking values in the space
$\mathcal{F}(V_\LLambda[0])$ of meromorphic
functions
of $\lambda\in\C$ with values in $V_\LLambda[0]$. 
Let $\delta_j$, $j=1,\dots,n$ be the standard basis of $\C^n$.
Then the
qKZB equations have the form
\[
v(\zz+p\delta_i)=K_i(\zz,\tau,p)v(\zz),
\qquad i=1,\dots, n.
\]
The qKZB operators
$K_i(\zz,\tau,p)$ act on the space
$\mathcal{F}(V_\LLambda[0])$
and are given by
\begin{eqnarray*}
K_j(\zz,\tau,p)&=&R_{j,j-1}(z_j-z_{j-1}+p,\tau)\cdots R_{j,1}(z_j-z_1+p,\tau)
\\ & &\Gamma_jR_{j,n}(z_j-z_n,\tau)\cdots R_{j,j+1}(z_j-z_{j+1},\tau).
\end{eqnarray*}
The operators $R_{j,k}(z,\tau)$ are defined by the formula
\[
R_{j,k}(z,\tau)\,v(\lambda)=R_{\Lambda_j,\Lambda_k}(z,\lambda-2\eta
\sum_{{l=1,l\neq j}}^{k-1}h^{(l)},\tau)\,v(\lambda),
\]
and $(\Gamma_jv)(\lambda)=v(\lambda-2\eta\mu)$ if $h^{(j)}v(\lambda)=
\mu v(\lambda)$ and is extended by linearity to 
$\mathcal{F}(V_\LLambda[0])$.

The qKZB system of difference equations
is compatible, i.e., we have
\begin{equation}\label{eq0}
K_j(\zz+p\delta_l,\tau,p)K_l(\zz,\tau,p)=
K_l(\zz+p\delta_j,\tau,p)K_j(\zz,\tau,p),
\end{equation}
for all $j,l$, as a consequence of the dynamical Yang--Baxter equations
satisfied by the $R$-matrices.
 We also consider the ``mirror'' qKZB operators
\begin{eqnarray*}
K_j^\vee(\zz,p,\tau)
&=&R^\vee_{j,j+1}(z_j-z_{j+1}+\tau,p)\cdots R^\vee_{j,n}(z_j-z_n+\tau,p)
\\ & &\Gamma_j
R^\vee_{j,1}(z_j-z_1,p)\cdots R^\vee_{j,j-1}(z_j-z_{j-1},p),
\end{eqnarray*}
with
\[
R^\vee_{j,k}(z,p)\,v(\lambda)=R_{\Lambda_j,\Lambda_k}(z,\lambda-2\eta
\sum_{{l=k+1,l\neq j}}^{n}h^{(l)},p)\,v(\lambda),
\]
The corresponding system of qKZB equations
\[v(\zz+\tau\delta_j)=K^\vee_j(\zz,p,\tau)\,v(\zz),\qquad j=1,\dots n,
\]
is also compatible. In fact, if 
we write $\vec x^\vee=(x_n,\dots,x_1)$ for any 
$\vec x=(x_1,\dots,x_n)\in \C^n$
and let $P:V_\LLambda\to V_{\LLambda^\vee}$
be the linear map sending $v_1\otimes\cdots\otimes v_n$ to $v_n\otimes
\cdots\otimes v_1$,
then we have, adding the dependence on $\LLambda$ in the notation,
\[
K^\vee_i(\zz,p,\tau;\LLambda)=
P^{-1}K_{n+1-i}(\zz^\vee,p,\tau;\LLambda^\vee)P.
\]

\subsection{Hypergeometric solutions}\label{suse-12}
In \cite{FTV2} we constructed a ``universal hypergeometric 
function'', which is a projective solution  of
the qKZB equations: it is 
 a function $u(\zz,\lambda,\mu,\tau,p)$, defined for generic values
of the parameters $\eta,\LLambda$, taking
values in $V_{\LLambda}[0]\otimes V_{\LLambda}[0]$ and obeying the
equations
\begin{eqnarray}\label{eq1}
u(\zz+\delta_jp,\tau,p)&=&K_j(\zz,\tau,p)\otimes D_j\, u(\zz,\tau,p),\notag\\
u(\zz+\delta_j\tau,\tau,p)&=&D^\vee_j\otimes K^\vee_j(\zz,p,\tau)\, u(\zz,\tau,p),\\
u(\zz+\delta_j,\tau,p)&=&u(\zz,\tau,p).
\notag
\end{eqnarray}
Here we view $u$ as taking values in the space of functions of
$\lambda$ and $\mu$ with values in  $V_{\LLambda}[0]\otimes V_{\LLambda}[0]$.
$K_j$ acts on the variable $\lambda$ and $K_j^\vee$ on the variable
$\mu$. The operators $D_j$, $D_j^\vee$ act by multiplication by
diagonal matrices $D_j(\mu)$, $D_j^\vee(\lambda)$, respectively.
For our purpose, the most convenient description of these matrices
is in terms of the function
\[
\alpha(\lambda)=\exp(-{\pi i\lambda^2/4\eta}).
\]
We have, for $j=1,\dots, n$,
\begin{eqnarray*}
D_j(\mu)&=&
\frac
{\alpha(\mu-2\eta(h^{(j+1)}+\cdots+h^{(n)}))}
{\alpha(\mu-2\eta(h^{(j)}+\cdots+h^{(n)}))}
\, e^{\pi i\eta\Lambda_j(\sum_{l=1}^{j-1}\Lambda_l-
                \sum_{l=j+1}^{n}\Lambda_l)},
\\
D^\vee_j(\lambda)&=&
\frac
{\alpha(\lambda-2\eta(h^{(1)}+\cdots+h^{(j-1)}))}
{\alpha(\lambda-2\eta(h^{(1)}+\cdots+h^{(j)}))}
\, e^{-\pi i\eta\Lambda_j(\sum_{l=1}^{j-1}\Lambda_l-
                \sum_{l=j+1}^{n}\Lambda_l)}.
\end{eqnarray*}
These operators are diagonal in the basis of $V_{\LLambda}[0]$
formed
by tensor products $e_I=e_{i_1}\otimes\cdots\otimes e_{i_n}$
of basis vectors of the $V_{\Lambda_k}$ so that $\sum (\Lambda_k-2i_k)=0$.

{}From $u$ one can construct projective solutions (eigenfunctions 
of the corresponding difference operators) by taking coefficients
of the basis vectors $e_I$. If 
$D_i(\mu)e_I=d_{i,I}(\mu)e_I$, $d_{i,I}(\mu)\in\C$, and
$u=\sum u^I\otimes e_I$ then for any fixed $I$ and $\mu$, the
function $\tilde v(\zz,\lambda)= u^I(\zz,\lambda,\mu,\tau,p)$ obeys
$\tilde v(\zz+p\delta_i)=d_{i,I}(\mu)K_i(\zz,\tau,p)\tilde v(\zz)$. It
follows that 
\[v(\zz,\lambda)=\prod_{i=1}^nd_{i,I}(\mu)^{-z_i/p}\tilde v(\zz,\lambda),
\] is a true solution
to the qKZB equations. The parameters $I$ and $\mu$ determine
the multipliers, as is easily seen from the explicit expression
for $u$ below:
\begin{equation}\label{multip}
v(\zz+\delta_i,\lambda)=d_{i,I}(\mu)^{-1/p}v(\zz,\lambda),\qquad
v(\zz,\lambda+1)= e^{-\frac{\pi i(\mu+2\eta m)}{2\eta}}v(\zz,\lambda).
\end{equation}
 The second
system of equations in \Ref{eq1} gives the 
monodromy of these solutions, see  \cite{FTV2}.

The explicit expression for $u$ is given by the following formulas.
\begin{eqnarray}\label{eq-expl}
u(\zz,\lambda,\mu,\tau,p)&=&
e^{-\frac{\pi i\lambda\mu}{2\eta}}
\int
\prod_{i,k}\Omega_{\eta\Lambda_k}(t_i-z_k,\tau,p)
\prod_{i<j}\Omega_{-2\eta}(t_i-t_j,\tau,p)\\
 & &
\sum_{I,J}\omega_I(t,\zz,\lambda,\tau)\omega^\vee_J(t,\zz,\mu,p)
dt_1\cdots dt_m
e_I\otimes e_J.\nonumber
\end{eqnarray}
The phase function $\Omega$ has the product formula
\[
\Omega_a(z,\tau,p)=\prod_{j,k=0}^\infty
\frac
{(1-e^{2\pi i(z-a+j\tau+kp)})(1-e^{2\pi i(-z-a+(j+1)\tau+(k+1)p)})}
{(1-e^{2\pi i(z+a+j\tau+kp)})(1-e^{2\pi i(-z+a+(j+1)\tau+(k+1)p)})}
\,.
\]
It is symmetric under exchanging $\tau$ and $p$ 
and obeys the functional equation 
\begin{equation}\label{e-feO}
\Omega_a(z+p,\tau,p)=e^{2\pi ia}\frac{\theta(z+a,\tau)}
{\theta(z-a,\tau)}\,\Omega_a(z,\tau,p).
\end{equation}
The weight functions $\omega_I$ are given by
\begin{gather*}
\omega_{(i_1,\dots,i_n)}(t_1,\dots,t_m,\zz,\lambda,\tau)=
\prod_{i<j}
\frac{\theta(t_i-t_j,\tau)}
{\theta(t_i-t_j+2\eta,\tau)}
\sum_{I_1,\dots,I_n}
\prod_{l=1}^n
\prod_{i\in I_l}
\prod_{k=1}^{l-1}
\frac{\theta(t_i-z_k+\eta\Lambda_k,\tau)}{\theta(t_i-z_k-\eta\Lambda_k,\tau)}
\\
\times\prod_{k<l}\prod_{i\in I_k,j\in I_l}
\frac{\theta(t_i-t_j+2\eta,\tau)}{\theta(t_i-t_j,\tau)}
\prod_{k=1}^{n}\prod_{j\in I_k}
\frac{\theta(\lambda\!+\!t_j\!-\!z_k\!-\!\eta\Lambda_k\!+\!2\eta
i_k\!-\!2\eta\sum_{l=1}^{k-1}(\Lambda_l\!-\!2i_l),\tau)}
{\theta(t_j-z_k-\eta\Lambda_k,\tau)}\,.
\end{gather*}
The summation is over all $n$-tuples $I_1,\dots,I_n$ of disjoint subsets of
$\{1,\dots,m\}$ such that $I_k$ has $i_k$ elements, $1\leq k\leq n$.
The theta function is here the first Jacobi theta function
\[
\theta(t,\tau)=-\sum_{j\in\Z}
e^{\pi i(j+\half)^2\tau+2\pi i(j+\half)(t+\half)}.
\]
The ``mirror'' weight functions are related to the weight functions
with the reversed order of factor. Indicating the dependence on the
highest weights explicitly, we have
$\omega_I^\vee(t,\zz,\mu,p,\LLambda)
=\omega_{I^\vee}(t,\zz^\vee,\mu,p,\LLambda^\vee),
$
where, as above, $(x_1,\dots,x_n)^\vee=(x_n,\dots,x_1)$.

The integral over $t_1,\dots,t_m$ of the 1-periodic integrand
in \Ref{eq-expl} is defined by analytic continuation from
a region of the space of parameter where the rule $\sum\Lambda_i=
2m$ does not hold: one starts
from the region $\mathrm{Im}(\eta)<0$, $\mathrm{Im}(\tau),
\mathrm{Im}(p),\mathrm{Im}(\eta\Lambda_i)>0$ for which the
integral is over the torus $(\R/\Z)^m$ and defines the integral
in general by analytic continuation.

\subsection{Remark}
In \cite{FTV2} we used only weight functions and no mirror
weight functions. Then the qKZB equations \Ref{eq1} only involve
qKZB operator and no mirror qKZB operators. The choices of this
paper make the qKZB heat equation more transparent.
The proof that $u$ obeys the relations \Ref{eq1} is the same
as the proof of Theorem 31 in \cite{FTV2}. Note however that
the conventions in the definitions of $D_j$ are different there.

\section{The qKZB heat equation}\label{se-2}
In this section we define a q-analogue of the KZB heat equation,
prove that it is compatible with the other qKZB equations and 
show that the differential KZB heat equation arises in the
semiclassical limit in the simplest non-trivial case.
 
The qKZB heat equation is an integral equation. The integration 
kernel is a contraction with the fundamental hypergeometric
solution.  The contraction is defined using 
the elliptic Shapovalov form.

\subsection{The elliptic Shapovalov form}\label{suse-21}
For $j=1,\dots,n$, $\mu,\tau\in\C$, $\mathrm{Im}\,\tau>0$, let
$Q^{\Lambda_j}(\mu,\tau):V_{\Lambda_j}\otimes V_{\Lambda_j}
\to\C$ be the bilinear form on $V_{\Lambda_j}$
with matrix elements $Q^{\Lambda_j}(\mu,\tau)(e_k\otimes e_l)=
\delta_{k,l}Q^{\Lambda_j}_k(\mu,\tau)$,
\begin{equation}\label{e-SH}
Q_k^{\Lambda_j}(\mu,\tau)=
\left(\frac{\theta'(0,\tau)}{\theta(2\eta,\tau)}\right)^{k}
\prod_{l=1}^k
\frac
{\theta(2\eta(\Lambda_j+1-l),\tau)\theta(2\eta l,\tau)}
{\theta(\mu+2\eta(\Lambda_j+1-k-l),\tau)\theta(\mu-2\eta l,\tau)}\, .
\end{equation}
Out of these bilinear forms we define a bilinear form 
$Q(\mu,\tau):V_{\LLambda}\otimes V_{\LLambda}\to\C$ on the
tensor product $V_\LLambda$: 
\[
Q(\mu,\tau)=Q^{\Lambda_1}(\mu,\tau)^{(1,n+1)} 
Q^{\Lambda_2}(\mu+2\eta h^{(1)},\tau)^{(2,n+2)} 
\cdots
Q^{\Lambda_n}(\mu+2\eta\sum_{j=1}^{n-1}h^{(j)},\tau)^{(n,2n)},
\]
and a bilinear form on the space of
functions of $\lambda$ with values in $V_{\LLambda}[0]$: if $f$
and $g$ are holomorphic functions from $\C$ to $V_{\LLambda}[0]$,
we set (if the integral converges)
\[
Q_\tau(f\otimes g)=\int Q(\mu,\tau) f(\mu)\otimes g(-\mu)
\alpha(\mu)\, d\mu.
\]
The integration is on the path $t\mapsto 2\eta t+\epsilon$, 
$-\infty< t<\infty$. This bilinear form is called the
elliptic Shapovalov form.

The main property of the elliptic Shapovalov form is that the $R$-matrix
is in a certain sense symmetric with respect to it, see Lemma \ref{lemma1.5}.
In particular it is a sort of contravariant form for the action of the
elliptic quantum group, see eq.~\Ref{firstclaim} below.

\subsection{Notation}\label{suse-notation}
To write the following formulae in the most transparent form
we will use the following notational conventions. 
Let for $k\in \Z_{\geq 1}$ and a complex vector space $V$,
 $\mathcal{F}_k(V)$ denote the
space of meromorphic functions of $k$ complex variables
with values in $V^{\otimes k}$.  For example $u(\zz,\tau,p)\in \mathcal{F}_2(V_\LLambda[0])$.
We also set $\mathcal{F}_0(V)=\C$ and $\mathcal{F}(V)=\mathcal{F}_1(V)$.
If $f\in\mathcal{F}_j(V)$ and $g\in \mathcal{F}_k(V)$, we define
$f\otimes g\in\mathcal{F}_{j+k}$ by 
$
(f\otimes g)(\lambda_1,\dots,\lambda_{j+k})
=f(\lambda_1,\dots,\lambda_j)\otimes g(\lambda_{j+1},\dots,\lambda_{j+k}).
$
If $A_i\in\mathrm{End}(\mathcal{F}(V))$ are difference operators
with meromorphic coefficients we write $A_1\otimes\cdots\otimes A_r\in
\mathrm{End}(\mathcal{F}_r(V))$
to denote the composition of the operators $A_i$, each acting on the $i$th
variable and the $i$th factor. We also use this notation if one of the
$A_i$ is an integral operator $Q:D\subset\mathcal{F}_2(V)\to\C$ 
of the form $f\mapsto\int Q(\mu)f(\mu,-\mu)d\mu$, $Q(\mu)\in (V\otimes V)^*$,
defined on
some subset $D$. Then $A_1\otimes\cdots\otimes A_r$ is defined on
some subset of $\mathcal{F}_r(V)$ and maps to $\mathcal{F}_{r-2}(V)$.

\subsection{The qKZB heat equation}\label{suse-22}
Let $T(\zz,\tau,p)$ be 
the integral operator on $V_{\LLambda}[0]$-valued
functions of one complex variable $\lambda$
\begin{equation}\label{e-he1}
T(\zz,\tau,p)v=(\alpha\otimes Q_{\tau+p}) u(\zz,\tau,\tau+p)\otimes v,
\end{equation}
where $\alpha$ is the operator of multiplication by the function
$\alpha(\lambda)$. More explicitly, if $\{e_I\}$ is a basis of 
$V_{\LLambda}[0]$ consisting of tensor products of basis vectors
and $u=\sum_{I,J}u_{I,J}e_I\otimes e_J$, $v=\sum v_Ie_I$,
$Q(\mu,\tau)(e_I\otimes e_J)=Q_{I}(\mu,\tau)\delta_{I,J}$, we have
\begin{equation}\label{e-he2}
T(\zz,\tau,p)v(\lambda)=\sum_I\left(\alpha(\lambda)\sum_{J}
\int u_{I,J}(\zz,\lambda,\mu,\tau,\tau+p)
Q_{J}(\mu,\tau+p)v_J(-\mu)\alpha(\mu)d\mu
\right)\,e_I.
\end{equation}

\medskip

\begin{thm}\label{t-1}
The equations 
\begin{eqnarray}\label{eq-ak}
 v(\zz+p\delta_j,\tau)&=&K_j(\zz,\tau,p)v(\zz,\tau),\qquad j=1,\dots, n,
\nonumber
\\ & & \\
 v(\zz,\tau)&=&T(\zz,\tau,p)v(\zz,\tau+p),\nonumber
\end{eqnarray}
are compatible, i.e., we have, in addition to \Ref{eq0},
\[
T(\zz+p\delta_j,\tau,p)K_j(\zz,\tau+p,p)
=K_j(\zz,\tau,p)T(\zz,\tau,p).
\]
\end{thm}

\medskip
The proof of this theorem is given in \ref{suse-24}.

A similar statement holds for the qKZB difference
operators $K^\vee_j$: they obey
the compatibility identities
\[
K^\vee_j(\zz+\tau\delta_l,p,\tau)K^\vee_l(\zz,p,\tau)=
K^\vee_l(\zz+\tau\delta_j,p,\tau)K^\vee_j(\zz,p,\tau),
\]
for all $j,l$, and the operator
\[
T^\vee(\zz,p,\tau)v=(Q_{\tau+p}\otimes\alpha)v\otimes u(\zz,\tau+p,p),
\]
obeys
\[
T^\vee(\zz+\tau\delta_j,p,\tau)
K^\vee_j(\zz,p+\tau,\tau)
=K^\vee_j(\zz,p,\tau)T^\vee(\zz,p,\tau),
\]
implying the compatibility of the mirror qKZB equations
\begin{eqnarray*}
 v(\zz+\tau\delta_j,p)&=&K^\vee_j(\zz,p,\tau)v(\zz,p),\qquad j=1,\dots, n,
\\
 v(\zz,p)&=&T^\vee(\zz,p,\tau)v(\zz,p+\tau).
\end{eqnarray*}

\medskip
\begin{corollary}\label{cor-1}
Let $U(\zz,\lambda,\mu,\tau,p)$ be the function
\[
U(\zz,\tau,p)=\alpha\otimes Q_{\tau+p}\otimes\alpha
\bigl(u(\zz,\tau,\tau+p)\otimes u(\zz,\tau+p,p)\bigr).
\]
Then $U$ 
obeys the system of equations \Ref{eq1}.
\end{corollary}

\medskip
\begin{conjecture}\label{conj}
Let $U(\zz,\tau,p)$ be the function defined in Corollary
\ref{cor-1}. Then
\[
U(\zz,\tau,p)=C u(\zz,\tau,p),
\]
 for some constant $C$.
\end{conjecture}

This conjecture is proved in the case where $\sum\Lambda_i=2$
(with $C=-e^{4\pi i\eta}/(2\pi\sqrt{4 i\eta})$).
This proof will be published elsewhere, \cite{FV3}.

Assuming the conjecture correct, we can obtain solutions
to the full system \Ref{eq-ak} as in \ref{suse-22}: for
arbitrary $I=(i_1,\dots,i_n)\in\Z_{\geq 0}^n$ with $\sum(\Lambda_k-2i_k)=0$
and $\mu\in\C$, let
\[
v(\zz,\lambda,\tau)=
e^{-\frac{i\pi \mu^2\tau}{4\eta p}}
\prod_{i=1}^nd_{i,I}(\mu)^{-z_i/p}
u^I(\zz,\lambda,\mu,\tau,p).
\]
Then the function $v(\zz,\tau):\lambda\mapsto v(\zz,\lambda,\tau)$
is a solutions of $\Ref{eq-ak}$. It also obeys \Ref{multip}
and 
$v(\zz,\lambda,\tau+1)
=e^{-\frac {i\pi\mu^2}{4\eta p}}
v(\zz,\lambda,\tau)$.

\medskip

\noindent{\bf Remark.} The Shapovalov pairing
$Q_\tau$ contains an integration $\int d\mu$ which we chose
to be the integral on the path $t\mapsto 2\eta t+\epsilon$.
This choice makes the integral convergent if $\mathrm{Im}(\eta)<0$
for a large class
of functions, thanks to the strong decay at infinity of the
Gaussian function $\alpha(\mu)$ on this path. This class contains
in particular our hypergeometric solutions.

It should be
however emphasized that the only properties of $\int d\mu$ that
are needed for this construction are that it be a linear form
on functions of $\mu$ invariant under translations by $2\eta$
times weights
of vectors in $V_{\Lambda_i}$, and that it be well defined on
a suitable class of functions. 

In particular, if the highest weights are rational with greatest common
denominator $d$, we
may replace the integral over $\mu$ by the sum over the set
$\{\lambda+2\eta k/d, k\in\Z\}$, so that the heat equation may be
viewed as a difference equation of infinite order.

\subsection{The case of integer highest weights}
If $\Lambda\in\Z_{\geq 0}$, let $S_\Lambda$ be the $\h$-submodule
of $V_\Lambda$ generated by $e_{\Lambda+1},e_{\Lambda+2},\dots$.
The $\Lambda+1$-dimensional quotient $V_\Lambda/S_\Lambda$ is
denoted by $L_{\Lambda}$. The classes $\bar e_0,\dots,\bar e_{\Lambda}$
of  $e_0,\dots,e_{\Lambda}$
build a basis of $L_{\Lambda}$. The space $L_{\Lambda}$ carries
 a one-dimensional
family of representations of the elliptic quantum group $E_{\tau,\eta}(sl_2)$,
see \cite{FV2}.

For integer highest weights $\Lambda_i,\Lambda_j$,
the $R$-matrix $R_{\Lambda_i,\Lambda_j}(z,\lambda,\tau)
$ preserves $S_{\Lambda_i}\otimes V_{\Lambda_j}$
and $V_{\Lambda_i}\otimes S_{\Lambda_j}$ \cite{FV2}, \cite{FTV1}.
Therefore it induces an endomorphism, 
still denoted $R_{\Lambda_i,\Lambda_j}(z,\lambda,\tau)$, of $L_{\Lambda_j}\otimes
L_{\Lambda_j}$. If $\Lambda_1,\dots,\Lambda_n\in\Z_{\geq 0}$
we then have a system of dynamical $R$-matrices and thus
a system of qKZB equations defined on $L_{\LLambda}[0]=
(L_{\Lambda_1}\otimes\cdots\otimes L_{\Lambda_n})[0]$.

A universal hypergeometric function $\hat u(\zz,\lambda,\mu,\tau,\eta)$
taking values in the tensor  product of 
finite dimensional modules $L_\LLambda[0]\otimes L_\LLambda[0]$ and obeying
\Ref{eq1} was found in \cite{MV}. It is defined by
$\hat u
(\zz,\lambda,\mu,\tau,\eta)=
\sum_{I,J}u_{I,J}(\zz,\lambda,\mu,\tau,p,\eta)\bar e_I\otimes \bar e_J$,
 where $u_{I,J}$ are the analytic continuation of the components of
the universal hypergeometric function 
for $V_{\vec \Lambda}[0] \otimes V_{\vec \Lambda}[0]$ which are shown to exist
for these values of $I,J$.

Then we can introduce a heat operator $\hat T(\zz,\tau,p)$ acting on 
functions with values in $L_\LLambda[0]$ by the same formula 
\Ref{e-he1}.
\medskip

\begin{thm}\label{t-1bis}
Suppose that $\Lambda_1,\dots,\Lambda_n$ are non-negative integers.
Then the equations 
\begin{eqnarray}\label{eq-ak1}
 v(\zz+p\delta_j,\tau)&=&K_j(\zz,\tau,p)v(\zz,\tau),\qquad j=1,\dots, n,
\nonumber
\\ & & \\
 v(\zz,\tau)&=&\hat T(\zz,\tau,p)v(\zz,\tau+p),\nonumber
\end{eqnarray}
for a function $v$ taking values in $L_\LLambda[0]$
are compatible, i.e., we have, in addition to \Ref{eq0},
\[
\hat T(\zz+p\delta_j,\tau,p)K_j(\zz,\tau+p,p)
=K_j(\zz,\tau,p)\hat T(\zz,\tau,p).
\]
\end{thm}

\medskip

The proof of this Theorem is contained in \ref{suse-76} below.

\subsection{Rational $\eta$}
A particularly interesting case is the case of 
integer highest weights and rational $\eta$.
 Let us for instance assume that $2N\eta=1$ for some
positive integer $N$ and suppose that the highest weights
$\Lambda_1,\dots,\Lambda_n$ are positive integers. 

If $N$ is large enough, then the qKZB operators may still be defined.
Indeed we have:
\medskip

\begin{lemma} 
Let $\Lambda_1,\Lambda_2$
be positive integers, and $N$ be a large enough integer.
Then the matrix elements of the $R$-matrix
$R_{\Lambda_1,\Lambda_2}
(z,\lambda,\tau;\eta)\in\mathrm{End}(L_{\Lambda_1}\otimes
L_{\Lambda_2})$ with respect to the basis
$\{\bar e_i\otimes \bar e_j\}$ are regular functions of $\eta$ at 
$2\eta=1/N$ for fixed generic values of $z,\lambda,\tau$.
\end{lemma}
\medskip

\noindent{\it Proof:} The $R$ matrix 
$R_{\Lambda_1,\Lambda_2}(z_1-z_2,\lambda,\tau;\eta)$,
for generic $\eta$,
is uniquely determined
up to normalization
by an intertwining condition for tensor products of 
$E_{\tau,\eta}(sl_2)$-modules $L_{\Lambda_i}(z_i)$, see \cite{FV2}.
The $E_{\tau,\eta}(sl_2)$ module $L_{\Lambda}(z)$ may be realized
for integer
$\Lambda$ as a symmetric tensor product of two-dimensional modules
\cite{FV2}, so that the
 matrix elements of $R_{\Lambda_i,\Lambda_j}$,
for a basis consisting of symmetrized tensor products of
basis vectors, can be expressed as polynomials in the matrix elements
of $R_{1,1}$. The latter 
matrix elements are known explicitly and are regular as
  functions of $\eta$. For $i=1,\dots,\Lambda$,
the basis vector
$\bar e_i$ of $L_\Lambda$
is proportional to  the symmetrized tensor products
of basis vectors of the two dimensional modules. The proportionality
constant is an elliptic factorial $\prod_{j=1}^i \theta(2\eta j)/\theta(2\eta)$
which is regular and non-zero at $2\eta N=1$ as long as $N>\lambda$.

Thus if $N>\mathrm{max}(\Lambda_1,\dots,\Lambda_n)$, the matrix elements
of the $R$-matrix are regular at $2\eta N=1$. $\square$

\medskip

Fix some generic complex number
$\epsilon$. Then we may consider the qKZB equations as equations for 
functions $v(\zz,\lambda)$, where the dynamical parameter $\lambda$
runs over the finite set $\{k/N+\epsilon\,|\,k\in \Z/2N\Z\}$.
Indeed, the coefficients of the qKZB operators are $1$-periodic functions
of $\lambda$ and the shifts of $\lambda$ in the difference operators
$\Gamma_i$ are integer multiples of $2\eta=1/N$. The shift by
the generic number $\epsilon$ ensures that on this finite set
no poles of the qKZB operators are encountered.  Thus, we get:

\medskip

\begin{proposition}\label{p-uno}
 Suppose that $N=(2\eta)^{-1}$ and
$\Lambda_1,\dots,\Lambda_n$ are positive 
integers, with $N$ large enough. Fix a generic complex number $\epsilon$.
  Let $F_N(\epsilon)$ be the space of functions 
$f:\frac1N\Z\to V_\LLambda[0]$ so that $f(\lambda+2)=f(\lambda)$.
Then the qKZB operators $K_i(\zz,\tau,p)$, $K^\vee_i(\zz,p,\tau)$
are well-defined endomorphisms of $F_N(\epsilon)$.
\end{proposition}

\medskip

In this situation we thus have a truly holonomic system, i.e.,
a compatible system of difference equations for functions
taking values in a finite dimensional vector space $F_N(\epsilon)$.
In order to define the heat equation we have to worry about the
fact that the universal hypergeometric function $\hat u$ 
is not defined for all values of the parameters. Recall that
$\hat u(\zz,\lambda,\mu,\tau,p)$ is also a meromorphic function of 
$\eta$. Let us say that 
$\eta$ is a regular point for $\hat u$ if $\hat u$ is regular
at this point for all $\lambda,\mu$ and all generic $\zz,\tau,p$.
\medskip

\begin{thm}\label{t-2}
Let $\eta$, $\LLambda=(\Lambda_1,\dots,\Lambda_n)$, $\epsilon$
 be as in Proposition \ref{p-uno} and assume that $\eta$ is
a regular point for $\hat u$. Then the heat operator
\[
\hat T_N(\zz,\tau,p)v(\lambda)=e^{-\frac{Ni\pi\lambda^2}2}
\sum_{k=0}^{2N-1} 
(1\otimes Q(\mu_k,\tau+p))
\hat u(\zz,\lambda,\mu_k,\tau,\tau+p)
\otimes v(-\mu_k)e^{-\frac{Ni\pi\mu_k^2}2},
\]
with $\mu_k=-\epsilon+k/N$, maps $F_N(\epsilon)$ to itself.
Moreover,
the equations for $v(\zz,\tau)\in F_N(\epsilon)$
\begin{eqnarray*}
 v(\zz+p\delta_j,\tau)&=&K_j(\zz,\tau,p)v(\zz,\tau),\qquad j=1,\dots, n,
\\
 v(\zz,\tau)&=&\hat T_N(\zz,\tau,p)v(\zz,\tau+p),\nonumber
\end{eqnarray*}
are compatible, i.e., we have, in addition to \Ref{eq0},
\[
\hat T_N(\zz+p\delta_j,\tau,p)K_j(\zz,\tau+p,p)
=K_j(\zz,\tau,p)\hat T_N(\zz,\tau,p),
\]
on $F_N(\epsilon)$.
\end{thm}

\medskip

The proof of this Theorem is contained in \ref{suse-76} below.

The description of the set of regular points for $\hat u$ will be
studied elsewhere. Here we only remark that in the case $n=1$,
$\Lambda_1=2$,
the point $\eta=1/2N$ is a regular point for
all $N\geq 3$, as can easily be checked since $\hat u$ is given by
a one-dimensional integral. In this case, Conjecture \ref{conj}
holds, see \cite{FV3}, namely, we have
\begin{eqnarray*}
\lefteqn{
\hat u(\zz,\lambda,\nu,\tau,p)=C_N
e^{-\frac{Ni\pi(\lambda^2+\nu^2)}2}
}
\\
&&\times
\sum_{k=0}^{2N-1}
(1\otimes Q(\mu_k,\tau+p)\otimes 1)
(\hat u(\zz,\lambda,\mu_k,\tau,\tau+p)
\otimes
\hat u(\zz,-\mu_k,\nu,\tau+p,p))
e^{-\frac{Ni\pi\mu_k^2}2},
\end{eqnarray*}
for all $\lambda\in\epsilon+\frac1N\Z$, $\nu\in\frac1N\Z$.
The constant is $C_N=\frac {ie^{2\pi i/N}}{S(N)}$, with
the Gauss sum
\[
S(N)=\sum_{k=0}^{2N-1}e^{-\frac{\pi i k^2}{2N}}
=(1-i)\sqrt N.
\]

\subsection{Proof of Theorem \ref{t-1}}\label{suse-24}
The proof is based on some identities involving $R$-matrices, $Q$ and $D_j$.
As above, we set $\alpha(\lambda)=\exp(-{\pi i\lambda^2/ 4\eta})$

\medskip
\begin{lemma}\label{lemma1.4}For any $\Lambda,\Mu$,
\[
\frac
{\alpha(\lambda-2\eta(h^{(1)}+h^{(2)}))}
{\alpha(\lambda-2\eta h^{(2)})}
R_{\Lambda,\Mu}(z+\tau,\lambda,\tau)
=
e^{-2\pi i\eta\Lambda,\Mu}
R_{\Lambda,\Mu}(z,\lambda,\tau)
\frac
{\alpha(\lambda-2\eta h^{(1)})}
{\alpha(\lambda)}
\, .
\]
\end{lemma}

\medskip
\begin{prf} One way to prove this lemma is to use the
functional realization (see \cite{FTV1}): The matrix
 $R_{\Lambda,\Mu}$ relate two bases of the same space of functions.
The basis elements are products of ratios of theta functions
and have therefore well-behaved transformation properties
under shifts of $z$ by $\tau$. The computation is straightforward
and will not be reproduced here.
\end{prf}
\medskip
\begin{lemma}\label{lemma1.5}
Let $\Lambda,\Mu\in\C$ and $v,w\in V_\Lambda\otimes V_\Mu$. 
Then
\begin{eqnarray*}
\lefteqn{
\langle  
Q^{\Lambda}(\mu+2\eta h^{(2)},\tau)
\otimes 
Q^{\Mu}(\mu,\tau) \, v,R_{\Lambda,\Mu}(z,-\mu,\tau)\,w
\rangle=
}\\
 & & =
\langle  
Q^{\Lambda}(\mu,\tau)
\otimes 
Q^{\Mu}(\mu+2\eta h^{(1)},\tau) R_{\Lambda,\Mu}(z,\mu+2\eta (h^{(1)}
+h^{(2)}),\tau) v,w
\rangle.
\end{eqnarray*}
\end{lemma}

\medskip

\begin{prf}
We first prove a version of this identity for $L$-operators.  Let
$L(\zeta,\lambda)\in{\mathrm{End}}
(\C^{2}\otimes V_\Lambda(z))$ be the $L$-operator
of the evaluation Verma module $V_\Lambda(z)$. We claim that,
for any $v_1,v_2\in\C^2\otimes V_\Lambda(z)$,
\begin{eqnarray}\label{firstclaim}
\lefteqn{
\langle  
Q^{{1}}(\mu+2\eta h^{(2)},\tau)
\otimes 
Q^{\Lambda}(\mu,\tau) \, v_1,L(\zeta,-\mu,\tau)\,v_2
\rangle=
}\notag\\
 & & =
\langle  
Q^{{1}}(\mu,\tau)
\otimes 
Q^{\Lambda}(\mu+2\eta h^{(1)},\tau) L(\zeta,\mu+2\eta (h^{(1)}
+h^{(2)}),\tau) v_1,v_2
\rangle.
\end{eqnarray}
We have $Q^1_0(\mu,\tau)=1$ and $Q^1_1(\mu,\tau)=
{\theta(2\eta)}{\theta(\mu-2\eta)^{-1}\theta(\mu)^{-1}}$. Define
the matrix elements of $L$ by $L(\zeta,\mu)e_j\otimes v=
\sum_{k=0,1}e_k\otimes L^k_j(\zeta,\mu)v$. Then the claim
is equivalent to 
\begin{eqnarray*}\lefteqn{
Q^\Lambda_k(\mu,\tau)
\langle e_k,
L^0_0(\zeta,-\mu)
e_k\rangle}\\
 &=&
Q^\Lambda_k(\mu+2\eta,\tau)
\langle
L^0_0(\zeta,\mu+2\eta(\Lambda-2k+1))
e_k,e_k\rangle,
\\
\lefteqn{Q^\Lambda_k(\mu,\tau)
\langle e_k,
L^0_1(\zeta,-\mu)
e_{k-1}\rangle}\\
 &=&
Q^1_1(\mu,\tau)
Q^\Lambda_{k-1}(\mu-2\eta,\tau)
\langle
L^1_0(\zeta,\mu+2\eta(\Lambda-2k+1))
e_k,e_{k-1}\rangle,
\\
\lefteqn{
Q^1_1(\mu+2\eta(\Lambda-2k),\tau)
Q^\Lambda_k(\mu,\tau)
\langle e_k,
L^1_0(\zeta,-\mu)
e_{k+1}\rangle}
\\
 &=&
Q^\Lambda_{k+1}(\mu+2\eta,\tau)
\langle
L^0_1(\zeta,\mu+2\eta(\Lambda-2k-1))
e_k,e_{k+1}\rangle,
\\
\lefteqn{Q^1_1(\mu+2\eta(\Lambda-2k),\tau)
Q^\Lambda_k(\mu,\tau)
\langle e_k,
L^1_1(\zeta,-\mu)
e_{k}\rangle}\\
 &=&
Q^1_1(\mu,\tau)
Q^\Lambda_{k}(\mu-2\eta,\tau)
\langle
L^1_1(\zeta,\mu+2\eta(\Lambda-2k-1))
e_k,e_{k}\rangle.
\end{eqnarray*}
These identities follow immediately from the explicit expressions
given in \cite{FV2} for the operators $L^j_k$ (called $a,b,c,d$ in
\cite{FV2}).

We now extend this result to the general case.
We use the intertwining property
of the $R$-matrix: let $\mathcal{L}_\Lambda, \mathcal{L}_\Mu$ be the
$L$-operators of $V_\Lambda(z_1)$, $V_\Mu(z_2)$, respectively.
Then%
\footnote{Here we omit the argument 
$\tau$ to 
shorten the notation} $R_{\Lambda,\Mu}(z_1-z_2,\mu)\in{\mathrm{End}}(V_\Lambda(z_1)
\otimes V_\Mu(z_2))$ is uniquely
determined up to a factor by the relation
\[
\mathcal{L}(\zeta,\mu)
R_{\Lambda,\Mu}(z_1-z_2,\mu-2\eta h^{(1)})^{(23)}
=
R_{\Lambda,\Mu}(z_1-z_2,\mu)^{(23)}
\mathcal{L}'(\zeta,\mu).
\]
The operators $\mathcal{L}$ and $\mathcal{L}'$ 
(giving the action of the quantum group
on the tensor product by using the coproduct and the opposite coproduct,
respectively)
are defined by
\begin{eqnarray*}
\mathcal{L}(\zeta,\mu)
 &=&
\mathcal{L}_\Lambda(\zeta,\mu-2\eta h^{(3)})^{(12)}
\mathcal{L}_\Mu(\zeta,\mu)^{(13)},\\
\mathcal{L}'(\zeta,\mu)
 &=&
\mathcal{L}_\Mu(\zeta,\mu-2\eta h^{(2)})^{(13)}
\mathcal{L}_\Lambda(\zeta,\mu)^{(12)}.
\end{eqnarray*}
The $R$-matrix normalized by the condition
$R_{\Lambda,\Mu}(z_1-z_2,\mu)e_0\otimes e_0=e_0\otimes e_0$.

In particular, if $v=w=e_0\otimes e_0$, the claim of the lemma
is correct for trivial reasons. We prove the general case by
induction: let us suppose that the lemma is proved for
$v,w$ of weight $\Lambda+\Mu-2j$, $j=0,\dots,k-1$, $k\geq1$. Now
it is known, see \cite{FV2}, that, for generic parameters,
the weight space $V_\Lambda(z_1)\otimes V_\Mu(z_2)[\Lambda+\Mu-2k]$
is spanned by vectors of the form
 $\mathcal{L}_1^0(\zeta,\lambda)\,x$
(or  $\mathcal{L}_1^0(\zeta,\lambda)\,x$), $\zeta\in\C$,
$x$ of weight $\Lambda+\Mu-2(k-1)$, and any
fixed generic $\lambda$. Indeed, if these vectors
did not span the weight space, they would be part of
a proper submodule, contradicting the irreducibility
of the tensor product.

By iterating \Ref{firstclaim}, we obtain
\begin{eqnarray*}
\lefteqn{\langle
Q^1(\mu+2\eta(h^{(2)}+h^{(3)}))\otimes
Q^\Lambda(\mu+2\eta h^{(3)})\otimes
Q^\Mu(\mu)
v_1,
\mathcal{L}(\zeta,-\mu)v_2\rangle}\\
 &=&\langle
Q^1(\mu)\otimes
Q^\Lambda(\mu+2\eta(h^{(1)}+ h^{(3)}))\otimes
Q^\Mu(\mu+2\eta h^{(1)})
\\
 & &\mathcal{L}'(\zeta,\mu+2\eta(h^{(1)}+h^{(2)}+h^{(3)}))
v_1,
v_2\rangle.
\end{eqnarray*}
In particular, if $v_1=e_0\otimes v$, $v_2=e_1\otimes w$,
one has
\begin{eqnarray*}
\lefteqn{
\langle
Q^{\Lambda}
(\mu+2\eta h^{(2)})\otimes
Q^\Mu(\mu)\,v
,
\mathcal{L}^0_1(\zeta,-\mu)\,w\rangle
}\\
 &=&
Q_1^1(\mu)
\langle
Q^\Lambda
(\mu+2\eta(
-1+h^{(2)}))
\otimes
Q^\Mu(\mu-2\eta)
\mathcal{L}'{}^{1}_0
(\zeta,\mu+2\eta(1+h^{(1)}+h^{(2)}))
\,v,
w\rangle.
\end{eqnarray*}
We turn to the proof of the induction step.
Suppose that $v,w$ have weight
$\Lambda+\Mu-2 k$, and write
$w=\mathcal{L}'{}^0_1(\zeta,\mu)x$.
Let us set $z=z_1-z_2$. Then
\begin{eqnarray*}
\lefteqn{
\langle
Q^\Lambda(\mu+2\eta h^{(2)})
\otimes
Q^\Mu(\mu)
\, v,
R_{\Lambda,\Mu}(z,-\mu)\,w\rangle
}
\\
&=&
\langle
Q^\Lambda(\mu+2\eta h^{(2)})
\otimes
Q^\Mu(\mu)
\, v,
R_{\Lambda,\Mu}(z,-\mu)\mathcal{L}'{}^0_1(\zeta,-\mu)
\,x\rangle
\\
&=&
\langle
Q^\Lambda(\mu+2\eta h^{(2)})
\otimes
Q^\Mu(\mu)
\, v,
\mathcal{L}^0_1(\zeta,-\mu)
R_{\Lambda,\Mu}(z,-\mu+2\eta)
\,x\rangle
\\
&=&
Q^1_1(\mu)
\langle
Q^\Lambda(\mu+2\eta (-1+h^{(2)}))
\otimes
Q^\Mu(\mu-2\eta)
\mathcal{L}'{}^1_0(\zeta,\mu+2\eta(1+h^{(1)}+h^{(2)}))
\, v,\\
 & &
R_{\Lambda,\Mu}(z,-\mu+2\eta)
\,x\rangle
\\
&=&
Q^1_1(\mu)
\langle
Q^\Lambda(\mu-2\eta)
\otimes
Q^\Mu(\mu+2\eta(-1+h^{(1)}))
\\
 & &R_{\Lambda,\Mu}(z,\mu+2\eta(-1+h^{(1)}+h^{(2)}))
\mathcal{L}'{}^1_0(\zeta,\mu+2\eta(1+\Lambda+\Mu
-2k))
\, v,
x\rangle.
\end{eqnarray*}
In the last step, we used the induction hypothesis.
The calculation continues by commuting $R$ with $\mathcal{L}'$,
and then by bringing $\mathcal{L}$ to the right. This last
part is similar to the above calculation read backwards, and will not be reproduced in detail.
One finally obtains, as desired,
\[
\langle
Q^\Lambda(\mu)\otimes Q^\Mu(\mu+2\eta h^{(1)})
R_{\Lambda,\Mu}
(z,\mu+2\eta(h^{(1)}+h^{(2)}))\,v,w\rangle.
\]
This completes the induction step and thus the
proof of the lemma.
\end{prf}
\medskip
\begin{lemma}\label{lemma1.6}
\begin{eqnarray*}
\alpha\, (D_j^\vee)^{-1}K_j(\zz,\tau,p+\tau)
&=&
K_j(\zz,\tau,p)\,\alpha\,
e^{-\pi i\eta\Lambda_j(\sum_{l\neq j}\Lambda_l)}
\\
\,\alpha\, D_j^{-1}K^\vee_j(\zz,p,\tau+p)
&=&
K^\vee_j(\zz,p,\tau)\,\alpha\,
e^{-\pi i\eta\Lambda_j(\sum_{l\neq j}\Lambda_l)}.
\end{eqnarray*}
\end{lemma}

\medskip

\begin{prf}
This is a straightforward consequence of the definition 
of the difference operators $K_j,K^\vee_j$ and of 
Lemma \ref{lemma1.4}.
\end{prf}
\medskip
\begin{lemma}\label{lemma1.7}
 Let $f$, $g$ be holomorphic functions from
$\C$ to $V_\Lambda$.
\begin{eqnarray*}
Q_{\tau+p}(f,K_j(\zz,\tau+p,p)g)
&=&
Q_{\tau+p}(D_j^{-1}K^\vee_j(\zz+p\delta_j,\tau+p,\tau)f,g)
e^{\pi i\eta\Lambda_j(\sum_{l\neq j}\Lambda_l)}
\\
Q_{\tau+p}(K^\vee_j(\zz,\tau+p,\tau)f,g)
&=&
Q_{\tau+p}(f,(D_j^\vee)^{-1}K_j(\zz+\tau\delta_j,\tau+p,p)g)
e^{\pi i\eta\Lambda_j(\sum_{l\neq j}\Lambda_l)}.
\end{eqnarray*}
\end{lemma}

\medskip

\begin{prf}
The proof of the first identity is given by using Lemma \ref{lemma1.5}
to bring the $R$-matrices in $K_j$ to the left, the translation invariance
of the integral to bring $\Gamma_j$ to the left, and \ref{lemma1.6}
to commute the resulting  $K_j^\vee(\zz+p\delta_j,\tau+p,-p)$ with 
$\alpha$. The proof of the second identity is similar.
\end{prf}
We can now complete the proof of Theorem \ref{t-1}.
Let $C_j=e^{\pi i\eta\Lambda_j(\sum_{l\neq j}\Lambda_l)}$ and
$v$ a function from $\C$ to $V_{\LLambda}[0]$.
\begin{eqnarray*}
\lefteqn{T(\zz+p\delta_j,\tau,p)K_j(\zz,\tau+p,p)\,v
=}\\ &=&
(\alpha\otimes Q_{\tau+p})u(\zz+p\delta_j,\tau,p+\tau)
\otimes
K_j(\zz,\tau+p,p)\,v
\\
&=&
C_j(\alpha\otimes Q_{\tau+p})(1\otimes D_j^{-1}
K^\vee_j(\zz+p\delta_j,\tau+p,\tau)
u(\zz+p\delta_j,\tau,p+\tau))\otimes v
\\
&=&C_j(\alpha\otimes Q_{\tau+p})((D_j^\vee)^{-1}
\otimes D_j^{-1}u(\zz+(p+\tau)\delta_j,\tau,\tau+p))\otimes v
\\
 &=& C_j(\alpha\otimes Q_{\tau+p})((D_j^\vee)^{-1}
K_j(\zz,\tau,\tau+p)\otimes 1\,u(\zz,\tau,\tau+p))\otimes v
\\
&=&
K_j(\zz,\tau,p)T(\zz,\tau,p)\,v.
\end{eqnarray*}

\subsection{Proof of Theorem \ref{t-1bis} and Theorem \ref{t-2}}
\label{suse-76}
Theorem \ref{t-1bis} can be proven in the same
way as Theorem \ref{t-1}.

There is however an apparent difficulty: the heat operator involves
the Shapovalov form which contains a sum over all components of
$u$, including those that a priori do not have a limit for integer
highest weights. The solution is provided by
 Theorem 2 of \cite{MV}: let us say that $I=(i_1,\dots,i_n)
\in\Z_{\geq 0}^n$ is admissible for $\LLambda$ if $i_a\leq \Lambda_a$
for all $a=1,\dots,n$. Then (a special case of) Theorem 2 states
that the components $u_{I,J}(\zz,\lambda,\mu,\tau,p)$
 such that
$I$ {\em or} $J$ is admissible, are regular functions of the highest
weights at $\LLambda$ for generic values of the other variables.
Moreover, $Q^{\Lambda}_k(\mu,\tau)$ vanishes if $k\geq \Lambda$,
cf.\ \eqref{e-SH}, so that the sum appearing in the compatibility
condition is effectively restricted to admissible indices.

Theorem \ref{t-2} is proven in the
same way as Theorem \ref{t-1} and \ref{t-1bis}.
In fact the only property of the integral over $\mu$ that is
used in the proof is the translation invariance. So the
same proof gives the compatibility relation in this case provided
the function of $\mu_k$ on the right-hand side is periodic
in $k$ with period $2N$. Now $u(\zz,\lambda,\mu,\tau,p)$ 
is $\exp(-\pi iN\lambda\mu)$
times a $2$-periodic function of $\lambda$ and $\mu$. So the
exponential factors combine into the expression
\[
e^{-\frac{i\pi N(\lambda+\mu_k)^2}{2}}=
e^{-\frac{i\pi N(\lambda-\epsilon+k/N)^2}{2}}.
\]
If $\lambda\in\epsilon +\frac1N\Z$, this expression is periodic
in $k$ with period $2N$.
The same argument shows that $T_N(\zz,\tau,p)v(\lambda)$ is
2-periodic in $\lambda$ for $\lambda\in\epsilon+\frac1N\Z$.
$\square$

\section{Semiclassical limit}\label{suse-23}
\begin{figure}
\rotatebox{-90}{\scalebox{.5}{\includegraphics{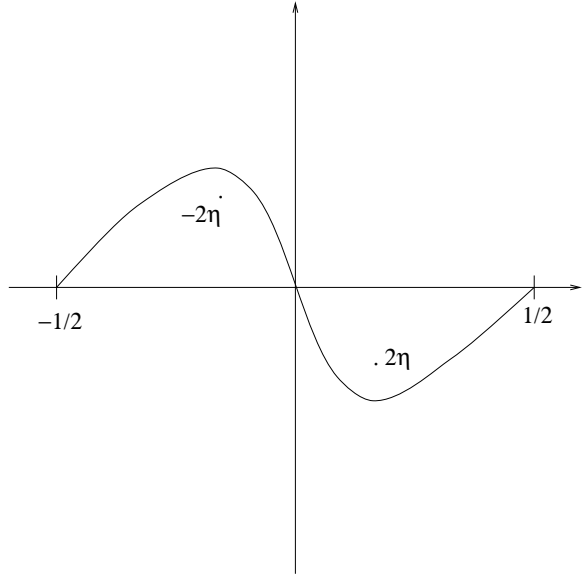}}}
\caption{The integration cycle $\gamma$. The points $\pm2\eta$ are
the singularities of the integrand}\label{Fig1}
\end{figure}

We consider here the semiclassical limit of our quantum heat equation
in the simplest non-trivial case and show that we do recover the KZB
heat equation in this limit. The case we consider is $n=1$, with
$\Lambda_1=2$. The qKZB equations for the dependence of $z_1$ are
trivial in this case and we can assume that $z_1=0$.  The zero weight
space is one-dimensional, and we identify it with $\C$ using the basis
$e_1$. Suppose that $v_\eta(\lambda,\tau)$ is a family of solutions of
the qKZB equations with parameters $\tau, p=-2\kappa\eta,\tau,\eta$,
parametrized by $\eta$ around zero. Assume that $v_\eta$ has an
asymptotic expansion $v_\eta(\lambda,\tau)=v_0(\lambda,\tau)+O(\eta)$
at
$\eta=0$.  We want to find the equations  satisfied by $v_0$. For this
we expand the qKZB heat equation
\begin{equation}\label{eq-qKZB1}
v_\eta(\lambda,\tau)=
\frac{-1}{4\pi\sqrt{i\eta}}
e^{-\frac{i\pi\lambda^2}{4\eta}}
\int u(\lambda,\mu,\tau,p,\eta)
\textstyle{\frac{\theta(4\eta,\tau+p)\theta'(0,\tau+p)}
{
 \theta(\mu+2\eta,\tau+p)
 \theta(\mu-2\eta,\tau+p)
}}
e^{-
\frac{i\pi\mu^2}
     {4\eta}
}
v_\eta(-\mu,\tau+p)\,d\mu,
\end{equation}
around $\eta=0$, setting $p=-2\kappa\eta$ and keeping $\tau,\kappa,\lambda$
fixed. The dependence of $\eta$  of the 
constant in front of the integral was chosen in such a way that
the semiclassical limit exists.
 The integration path is $t\mapsto \mu=\eta t$ ($t\in\R$).
The hypergeometric solution $u$ is independent of $z$ in this case and
is given by the formula:
\[
u(\lambda,\mu,\tau,p,\eta)=e^{-\frac{i\pi\lambda\mu}{2\eta}}
\int_\gamma
\Omega_{2\eta}(t,\tau,\tau+p)
\frac
{\theta(\lambda+t,\tau)}
{\theta(t-2\eta,\tau)}
\frac
{\theta(\mu+t,\tau+p)}
{\theta(t-2\eta,\tau+p)}
dt.
\]
The integration cycle $\gamma$ is depicted in Fig.~\ref{Fig1}.

\medskip
\begin{thm}
Suppose that $v_\eta(\lambda,\tau)$ is a family of solutions
of \Ref{eq-qKZB1} with an asymptotic expansion
$v_\eta(\lambda,\tau)=v_0(\lambda,\tau)+\eta v_1(\lambda,\tau)
+\cdots$, then
$v(\lambda,\tau)=v_0(\lambda,\tau)/\theta(\lambda,\tau)$
obeys the KZB heat equation
\[
2\pi i\kappa\frac
{\partial v}{\partial\tau}=
\frac
{\partial^2v}
{\partial\lambda^2}
-2\wp(\lambda,\tau)v
+c(\tau)v,
\]
for some $c(\tau)$ independent of $\lambda$.
\end{thm}\medskip

\noindent{\it Proof:} 
The integral on the right-hand side of \Ref{eq-qKZB1} has
the form 
\[
I_\eta=\frac{i}{\sqrt{4 i\eta}}\int _{-\eta\infty}^{\eta\infty}
e^{-\frac
{i\pi}
{4\eta}
(\lambda+\mu)^2
}
g(\lambda,-\mu,\eta)d\mu.
\]
This integral has the asymptotic expansion as $\eta\to 0$
\[
I_\eta=g(\lambda,\lambda,0)+
\eta\left(\frac1{i\pi}\frac{\partial^2}{\partial\mu^2}\bigg|_{\mu=\lambda}
g(\lambda,\mu,0)
+\frac{\partial}{\partial\eta}\bigg|_{\eta=0} 
g(\lambda,\lambda,\eta)\right)+O(\eta^2).
\]
To compute the various terms of this expression, we first notice that
the integration cycle in $u$ is pinched by the singularities as
$\eta\to 0$. The integral defining $u$ can then be expressed as a
divergent (as $\eta\to 0$) part given by $2\pi i$ times the residue at
$t=2\eta$ plus the integral on a cycle $\bar\gamma$ which stays away
from the singularities.

To compute the residue we introduce $\tilde \Omega$ by
\begin{eqnarray*}
\Omega_{2\eta}(t,\tau,\tau+p)&=&
\frac
{1-e^{2\pi i(t-2\eta)}}
{1-e^{2\pi i(t+2\eta)}}
\tilde \Omega_{2\eta}(t,\tau,\tau+p)
\\
&=&
(t-2\eta)
\frac{2\pi i}{e^{8\pi i\eta}-1}
\tilde\Omega_{2\eta}(2\eta,\tau,\tau+p) + O((t-2\eta)^2).
\end{eqnarray*}
As $\eta\to 0$, $\tilde\Omega_{2\eta}(2\eta,\tau,\tau-2\kappa\eta)$ is
regular and converges to $1$. 

We then have
\begin{eqnarray*}
g(\lambda,\mu,\eta)
&=&
\frac{2\pi i}
{e^{8\pi i\eta}-1}
\tilde\Omega_{2\eta}
(2\eta,\tau,\tau+p)
\frac{\theta(\lambda+2\eta,\tau+p)\theta(4\eta,\tau+p)
v_\eta(\mu,\tau+p)}
{\theta'(0,\tau)\theta(\mu+2\eta,\tau+p)}
\\
 & &
-\frac1{2\pi i}\int_{\bar\gamma}
\Omega_{2\eta}(t,\tau,\tau+p)
\frac
{\theta(\lambda+t,\tau)}
{\theta(t-2\eta,\tau)}
\frac
{\theta(\mu+t,\tau+p)}
{\theta(t-2\eta,\tau+p)}
dt\\
 & &\times
\frac{\theta(4\eta,\tau+p)\theta'(0,\tau+p)}
{
 \theta(\mu+2\eta,\tau+p)
 \theta(\mu-2\eta,\tau+p)
}v_\eta(\mu,\tau+p).
\end{eqnarray*}
{}From these formulae we can compute the various terms:
\[
g(\lambda,\lambda,0)=
\,v_0(\lambda,\tau),
\]
\[
\frac{\partial^2}{\partial^2\mu}\bigg|_{\mu=\lambda} g(\lambda,\mu,0)=
\theta(\lambda,\tau)\partial_\lambda^2
\left(
\frac {v_0(\lambda,\tau)}
{\theta(\lambda,\tau)}
\right).
\]
Finally
\begin{eqnarray*}
\frac{\partial}{\partial\eta}\bigg|_{\eta=0}g(\lambda,\lambda,\eta)
&=&C_1(\tau)\,v_0(\lambda,\tau)
+
\eta
\frac
{\partial}
{\partial\eta}
\bigg|_{\eta=0}
v_\eta(\lambda,\tau)
\\
 & &
-2\kappa\,\theta(\lambda,\tau)
\frac{\partial}{\partial\tau}
\left(
\frac {v_0(\lambda,\tau)}
{\theta(\lambda,\tau)}
\right)
\\
 & &-
\frac2{\pi i}
\int_{\bar\gamma}
\frac{\theta(t+\lambda,\tau)\theta(t-\lambda,\tau)}
{\theta(t,\tau)^2\theta(\lambda,\tau)^2}
dt
\,v_0(\lambda,\tau).
\end{eqnarray*}
Here $C_1(\tau)$ is some scalar function independent of $\lambda$.
Using the identity
\[
\frac{\theta(t+\lambda,\tau)\theta(t-\lambda,\tau)}
{\theta(t,\tau)^2\theta(\lambda,\tau)^2}
=\frac{1}{\theta'(0,\tau)^2}\left(\wp(\lambda,\tau)-
\wp(t,\tau)\right),
\]
we see that the right-hand side of \Ref{eq-qKZB1} is
\begin{eqnarray*}
\lefteqn{v_0(\lambda,\tau)
+
\eta\frac{\partial}{\partial\eta}\bigg|_{\eta=0}v_\eta(\lambda,\tau)}\\
 &&+\eta\theta(\lambda,\tau)\left(
\frac1{i\pi}\frac{\partial^2}{\partial\lambda^2}
-2\kappa\frac{\partial}{\partial\tau}-\frac2{\pi i}\wp(\lambda,\tau)
+c(\tau)\right)\frac
{v_0(\lambda,\tau)}
{\theta(\lambda,\tau)}+O(\eta^2),
\end{eqnarray*}
for some function $c(\tau)$ independent of $\lambda$.

Since the first two terms also appear on the left-hand side,
the proof is complete. $\square$

\section{Conformal blocks}\label{se-4}

In this section, we introduce,  in the simplest case of one marked point,
a difference analogue
the vector bundle of conformal blocks. We begin by reviewing the differential
case. The vector bundle of conformal blocks 
 is, in this case, a vector bundle on the moduli space $\mathcal{M}_{1,1}$
of genus one curves with one marked point.
The projectivization of this vector bundle carries a connection
given by the KZB differential operator.
We then give a difference analogue of this vector bundle. It has a 
 (discrete) connection,
which is now given by the qKZB heat operator $T$.

\subsection{The differential case}
Let $\mathfrak{g}$ be a simple complex Lie algebra with Cartan subalgebra
$\mathfrak{h}$ and root space decomposition
$\mathfrak{g}=\mathfrak{h}\oplus\oplus_{\alpha\in\Delta}
\mathfrak{g}_\alpha$. Let the non-degenerate invariant bilinear form
$(\ ,\ )$ on $\mathfrak g\simeq
\mathfrak g^*$ be normalized so that the highest root $\theta$ obeys
$(\theta,\theta)=2$. Let $\kappa$ be an integer larger than or equal to the
dual Coxeter number $h^\vee$
of $\mathfrak g$, and $\Lambda\in \mathfrak h^*$
be a dominant integral weight, 
so that $(\theta,\Lambda)\leq \kappa-h^\vee$. Denote by $L_\Lambda$ the
irreducible $\g$-module of highest weight $\Lambda$.

To these data one associates a holomorphic 
vector bundle of conformal blocks
on the moduli space ${\mathcal{M}}_{1,1}$ of genus one complex
curves with one marked point \cite{TUY}. Its projectivization carries
a canonical flat connection. The fiber over a point
may be defined as a space of
coinvariants for the Lie algebra of $\mathfrak g$-valued
rational functions on the curve whose poles are at the marked point,
acting on the irreducible affine Kac--Moody Lie algebra module
of highest weight $\Lambda$ and level $\kappa-h^\vee$.

An explicit description \cite{FW}
 of this bundle, which for our purposes can be
taken as a definition,
may be obtained by viewing ${\mathcal{M}}_{1,1}$ as the 
quotient of the upper half plane $H_+$ by $\mathrm{SL}(2,\Z)$.
We may then regard the 
 vector bundle $E_{\kappa,\Lambda}$ of conformal blocks as
an $\mathrm{SL}(2,\Z)$-%
equivariant vector bundle over $H_+$. Let $L_\Lambda[0]=
\{v\in L_\Lambda\,|\,\mathfrak{h}v=0\}$ be the zero weight space
space of $L_\Lambda$. It carries a natural linear action of the
Weyl group $W$ of $\mathfrak g$.
The fiber of $E_{\kappa,\Lambda}$
over $\tau\in H_+$ 
is then to the space of holomorphic maps
$v:\mathfrak h \to L_\Lambda[0]$ such that
\begin{enumerate}
\item[(i)] $v(\lambda+q_1+q_2\tau)=
\exp(-\pi i\kappa(q_2,q_2)\tau-2\pi i\kappa(q_2,\lambda))v(\lambda)$,
for all $\lambda\in\mathfrak h$ and $q_1,q_2$ in the coroot lattice
$Q^\vee$.
\item[(ii)] $v(w\cdot\lambda)=\epsilon(w)w\cdot v(\lambda)$ 
for all $w\in W$,
where $\epsilon:W\mapsto\{\pm1\}$ is the homomorphism sending
reflections to $-1$.
\item[(iii)] For all roots $\alpha$, 
$x\in\mathfrak g_\alpha$, and integers $l\geq0,r,s$,
the map 
$v$ obeys the vanishing condition 
\[
x^lv(\lambda)=O\bigl((\alpha(\lambda)-r-s\tau)^{l+1}\bigr),
\]
 as $\alpha(\lambda)\to r+s\tau$.
\end{enumerate}

The action of $\mathrm{SL}(2,\Z)$ on the base may be lifted
to an action on the bundle:
let 
$g=\left({{\begin{array}{cc}a&b\\ c&d\end{array}}}\right)
\in
\mathrm{SL}(2,\Z)$ act on $H_+$ by $\tau\mapsto g\cdot\tau=
(a\tau+b)/(c\tau+d)$. Then we have isomorphisms
$\psi_g(\tau):E_{\kappa,m}(\tau)\to E_{\kappa,m}(g\cdot\tau)$
given by
\[
\psi_g(\tau)v(\lambda)=
e^{\frac{\pi i\kappa}2 c(c\lambda+d)\lambda^2}
v((c\tau+d)\lambda),
\]
obeying the cocycle condition $\psi_{gh}(\tau)=\psi_g(h\cdot\tau)\psi_h(\tau)$.
Denote by $\pi:H_+\to\mathcal{M}_{1,1}$
the canonical projection.
Local holomorphic sections of the vector bundle of conformal blocks on
an open set $U\subset \mathcal{M}_{1,1}$ are then the same as
holomorphic sections $v$ of $E_{\kappa,m}$ on $\pi^{-1}(U)$
so that $v(g\cdot\tau)=\psi_g(\tau)^{-1}v(\tau)$. In other words,
they are holomorphic functions $v(\lambda,\tau)$ on $\C\times
p^{-1}(U)$ obeying (i)-(iii) for each fixed $\tau$ and such that
\[
v\left(\frac\lambda{c\tau+d},\frac{a\tau+b}{c\tau+d}\right) 
=e^{-\frac{\scriptstyle{\pi i\kappa c\lambda^2}}{\scriptstyle{2(c\tau+d)}}} 
v(\lambda,\tau).
\]
The projectivization of this 
vector bundle carries a holomorphic connection, 
and horizontal sections
may be constructed by an elliptic version of hypergeometric integrals
\cite{FV1}.

 We describe
here the  connection in the case of $sl(2,\C)$.
If $\mathfrak{g}=sl(2,\C)$ and $\Lambda=m\alpha$, $m=0,1,\dots$, 
then $L_\Lambda[0]$ is one dimensional. Let us chose a basis of
$L_\Lambda[0]$ and identify $\mathfrak h\simeq\mathfrak h^*$ 
with $\C$ via the basis
$\alpha/2$.  Then $E_{\kappa,\Lambda}(\tau)=E_{\kappa,2m}(\tau)$ 
consists of holomorphic 
functions $v(\lambda)$ on the complex plane so that 
(i) $v(\lambda+2r+2s\tau)=\exp(-2\pi i\kappa(s^2\tau +s\lambda))
v(\lambda)$,
(ii) $v(-\lambda)=(-1)^{m+1}v(\lambda)$, (iii) $v$ is divisible
 by $\theta(\lambda,\tau)^{m+1}$ in the ring of holomorphic functions.

If $\kappa\geq 2m+2$, we have 
$E_{\kappa,2m}(\tau)=\theta(\lambda,\tau)^{m+1}
\Theta_{\kappa-2m-2}(\tau)^W$,
where $\Theta_\kappa(\tau)^W$ is the $\kappa+1$-dimensional
space of holomorphic even functions
obeying (i). Otherwise $E_{\kappa,2m}(\tau)$ is trivial.

It follows that
\[
\mathrm{dim}(E_{\kappa,2m}(\tau))=\left\{
\begin{array} {rl}
\kappa-2m-1,& \mathrm{if}\;\kappa\geq 2m+2,\\
0,&\mathrm{otherwise.}
\end{array}\right.
\]
The connection on $E_{\kappa,2m}$
is defined by its covariant derivative $\Gamma(U,E_{\kappa,2m})
\to \Gamma(U,E_{\kappa,2m})\otimes\Omega^1(U)$ on local holomorphic
sections:
\[
\nabla v(\lambda,\tau)=\left(\partial_\tau-\frac1{2\pi i\kappa}
\bigl(\partial_\lambda^2-m(m+1)\wp(\lambda,\tau)\bigr)-\eta(\tau)^{-1}
\partial_\tau\eta(\tau)\right) v(\lambda,\tau) \,d\tau.
\]
Here $\wp$ is the Weierstrass elliptic function with periods 1 and
$\tau$ and 
\[
\eta(\tau)=e^{\pi i\tau/12}\prod_{j=1}^\infty(1-e^{2\pi i j\tau})
\]
is the Dedekind $\eta$-function.\footnote{The connection, being
on the projectivization, is really defined
up to adding a  multiple of the identity. We have chosen it here so that it
defines a connection on the vector bundle over $\mathcal{M}_{1,1}$}
In spite of the poles of the $\wp$ function, this connection is well-defined
on $E_{\kappa,2m}$ as can be seen by noticing that the poles cancel
in the expression of the 
induced connection $\theta^{-m-1}\circ\nabla\circ\theta^{m+1}$ 
on $\Theta_{\kappa-2m-2}^W$. The fact that $\nabla$ preserves (i)
and (ii) is easily checked.

The connection $\nabla$ is $\mathrm{SL}(2,\Z)$-equivariant, 
in the following sense: if $U\subset H_+$ is an 
$SL(2,\Z)$-invariant open set, and $g\in \mathrm{SL}(2,\Z)$, we have
the pull-back $g^*:\Gamma(U,E_{\kappa,2m})\to\Gamma(U,E_{\kappa,2m})$,
sending a section $v(\tau)$ to $\psi_g(\tau)^{-1}v(g\cdot\tau)$.
We may extend $g^*$ to $\Gamma(U,E_{\kappa,2m})\otimes\Omega^1(U)$ by
tensoring with the pull-back of differential forms. Then
$g^*\circ\nabla=\nabla\circ g^*$. Therefore the connection is well-defined
on the vector bundle of conformal blocks on $\mathcal{M}_{1,1}$.

\medskip

\noindent{\it Example.} If $m=0$, $\nabla$ is essentially the differential operator
of the heat equation. The theta functions
\[
\theta_{j,\kappa}(\lambda,\tau)=\sum_{r\in\Z+j/2\kappa}
e^{2\pi i\kappa (r^2\tau+r\lambda)},\qquad j\in\Z/2\kappa\Z,
\]
form a basis of $\Theta_{\kappa}(\tau)$ for fixed $\tau$, and obey the
heat equation $2\pi i\kappa\partial_\tau\theta_{j,\kappa}
=\partial_\lambda^2\theta_{j,\kappa}$. Moreover, we have $\theta_{j,\kappa}
(-\lambda,\tau)
=\theta_{-j,\kappa}(\lambda,\tau)$. It follows that the functions
\begin{equation}\label{e-s1}
v_j(\lambda,\tau)=\eta(\tau)^{-1}\bigl(\theta_{j+1,\kappa}(\lambda,\tau)
-\theta_{-j-1,\kappa}(\lambda,\tau)\bigr),\qquad j=0,1,\dots,\kappa-2,
\end{equation}
form a basis of the space of horizontal sections.

See \cite{FV1} for the case of arbitrary $m$.

\subsection{The difference case}
Let us turn to the difference case (for $sl(2,\C)$). 
We describe a difference analogue of $E_{\Lambda,2m}$,
 a holomorphic vector bundle
$E_{\Lambda,2m,\eta}$ on $H_+$ which is preserved by the qKZB heat operator.
We fix a generic $\eta$ in the lower half plane. Guided
by the semiclassical analysis of Sect.~\ref{suse-23},
we suppose that $-p/2\eta=\kappa$
is an integer $\geq2$ and consider the qKZB heat operator
 \eqref{e-he1} for $n=1$, $z_1=0$, $\Lambda_1=2m$.

We start with the somewhat trivial but instructive case $m=0$,
and write $T_{\kappa,0}(\tau)=T(z=0,\tau,p=-2\eta\kappa)$. Here
the qKZB heat operator is 
\[
T_{\kappa,0}(\tau)
v(\lambda)=\int_{2\eta\R} e^{-\frac{\pi i}{4\eta}(\lambda+\mu)^2}v(-\mu)\,d\mu.
\]
The integral is over the path $t\mapsto 2\eta t$, $-\infty<t<\infty$.

We define $E_{\kappa,2m=0,\eta}=E_{\kappa,0}$ to be the holomorphic vector bundle
of odd theta  functions, as in the differential case: the fiber over $\tau\in H_+$ is
$E_{\kappa,0,\eta}(\tau)=\{f\in\Theta_\kappa(\tau)\,|\, f(-\lambda)=-f(\lambda)\}$

\medskip

\begin{thm} Let $\kappa\geq 2$ and suppose that $\mathrm{Im}\,\eta<0$,
 $\mathrm{Im}\,\tau>0$.
Then $T_{\kappa,0}(\tau)$ maps $E_{\kappa,0,\eta}(\tau-2\eta\kappa)$ to
$E_{\kappa,0,\eta}(\tau)$.
\end{thm}

\medskip 

This theorem is based on the identity
\[
\theta_{j,\kappa}(\lambda,\tau)=\frac i{\sqrt{4i\eta}}
\int_{2\eta\R}e^{-\frac{i\pi}{4\eta}(\lambda+\mu)^2}
\theta_{j,\kappa}(-\mu,\tau-2\eta\kappa)\,d\mu,\qquad j\in\Z/2\kappa\Z,
\]
which gives the action of $T_{\kappa,0}(\tau)$ on the basis 
$\theta_j-\theta_{-j}$, $j=1,\dots,\kappa-1$,
of $\Theta_\kappa(\tau-2\eta\kappa)$.

\medskip

Let us now turn to the case of general $m$.
To compare with the classical limit we consider the qKZB operator
for the quotient $v$ of the dependent function by 
$\prod_{j=1}^m\theta(\lambda+2\eta j,\tau)$,
i.e., we set
\[
T_{\kappa,m}(\tau)=\phi_m(\tau)^{-1}\circ T(z=0,\tau,p=\tau-2\eta\kappa)
\circ \phi_m(\tau-2\eta\kappa),
\]
where $\phi_m(\tau)$ is the operator of multiplication by 
the function $\lambda\mapsto\prod_{j=1}^m\theta(\lambda+2\eta j,\tau)$.

\medskip

\noindent{\it Example.} If $m=1$, the qKZB operator for $v$ is
$v\mapsto T_{\kappa,1}(\tau)v$ is
\[
T_{\kappa,1}(\tau)v(\lambda)=\alpha(\lambda)
\int_{2\eta\R} V(\lambda,\mu,\tau,\tau-2\eta\kappa)\alpha(\mu)v(-\mu)d\mu,
\]
with kernel
\[
V(\lambda,\mu,\tau,\sigma)=c\,
e^{-\frac{\pi i\lambda\mu}{2\eta}}
\int_\gamma \Omega_{2\eta}(t,\tau,\sigma)
\frac{\theta(\lambda+t,\tau)\theta(\mu+t,\sigma)}
       {\theta(t-2\eta,\tau)\theta(\lambda+2\eta,\tau)
         \theta(t-2\eta,\sigma)\theta(\mu+2\eta,\sigma)}\,dt,
\]
for some $c=c(\tau,\sigma)$ independent of $\lambda,\mu$. The integration
cycle is depicted in Fig.~\ref{Fig1}.

\medskip

Let $E_{\kappa,2m,\eta}(\tau)$ be the space of holomorphic functions
so that 
 \begin{enumerate}
\item[(i)] $v(\lambda+2r+2s\tau)=
\exp(4\pi i\eta m(m+1) s-2\pi i\kappa (s^2\tau +s\lambda))
v(\lambda)$,
\item[(ii)]
 $v(-\lambda)=(-1)^{m+1}\prod_{j=1}^m\frac{\theta(\lambda+2\eta j,\tau)}
{\theta(\lambda-2\eta j,\tau)}v(\lambda)$, 
\item[(iii)] $v$ is divisible
 by $\prod_{j=0}^m\theta(\lambda-2\eta j,\tau)$ 
in the ring of holomorphic functions.
\end{enumerate}

Alternatively (and more simply), 
$E_{\kappa,2m,\eta}(\tau)$ is the space of functions of the
form $\prod_{j=0}^m
\theta(\lambda-2\eta j,\tau)\,\varphi(\lambda)$, with 
$\varphi\in\Theta_{\kappa-2m-2}(\tau)^W$.
In particular, $E_{\kappa,2m,\eta}(\tau)$ has the same
 dimension as the space $E_{\kappa,2m}(\tau)$ appearing in
the differential case. Let $E_{\kappa,2m,\eta}=\cup_{\tau\in H_+}
 E_{\kappa,2m,\eta}(\tau)$.
It is naturally a holomorphic vector bundle over $H_+$.

\begin{thm} Let $m,\kappa\in\Z_{\geq0}$,
$\kappa\geq 2m+2$ and suppose that $\mathrm{Im}\,\eta<0$,
 $\mathrm{Im}\,\tau>0$. Then
$T_{\kappa,m}(\tau)$ maps $E_{\kappa,2m,\eta}(\tau-2\eta\kappa)$
to $E_{\kappa,2m,\eta}(\tau)$.
\end{thm}

\medskip

\noindent{\it Proof:} This theorem is a corollary of the results of
\cite{FV4}.
We give here the proof in the simplest case $m=1$. The proof of the general
case is similar.
Let $v\in E_{\kappa,2,\eta}(\tau-2\eta\kappa)$,
and set $\tilde v=T_{\kappa,1}(\tau)v$. Properties (i), (ii) for $\tilde v$ can be checked
straightforwardly, by using the identities 
\[
\theta(\lambda+2,\tau)=\theta(\lambda,\tau),
\qquad
\theta(\lambda+2\tau,\tau)=e^{-4\pi i(\lambda+\tau)}\theta(\lambda,\tau),
\qquad
\theta(-\lambda,\tau)=-\theta(\lambda,\tau),
\]
obeyed by $\theta$ and translating the integration variable in the
integral over $\mu$. The latter involves moving the integration contour,
which presents no problem as the vanishing condition (iii) for $v$ 
guarantees that the integrand has no poles.
Let us check that $\tilde v$ is holomorphic and
obeys (iii). As the zeros of $\theta(\lambda,\tau)$ 
are simple and on the lattice
$\Z+\tau\Z$, $\tilde v$ is regular except possibly for
simple poles at $-2\eta+\Z+\tau\Z$. 

We claim that 
$\tilde v$ vanishes at $\lambda=r+s\tau$ and at $\lambda=2\eta+r+s\tau$
for all $r,s\in\Z$. Then (ii) implies that $\tilde v$ is regular
at the points $-2\eta+\Z+\tau\Z$ (and thus everywhere), 
and that $\tilde v$ is divisible
by $\theta(\lambda,\tau)\theta(\lambda-2\eta,\tau)$.

Since $\tilde v$ obeys (i), it is sufficient 
to prove the claim for $r,s\in\{0,1\}$

It follows from (ii) that $\tilde v(0)=0$ and, in
conjunction with (i), also $\tilde v(r+s\tau)=0$, $r,s\in\{0,\pm 1\}$.
For example, we have
\[ \tilde v(-\tau)=\frac{\theta(\tau+2\eta,\tau)}{\theta(\tau-2\eta,\tau)}
\,\tilde v(\tau)=
\frac{e^{-8\pi i\eta}\theta(-\tau+2\eta,\tau)}{\theta(\tau-2\eta,\tau)}
\,\tilde v(\tau)=-
e^{-8\pi\eta}\tilde v(\tau).
\]
On the other hand, (i) implies $\tilde v(-\tau)=e^{-8\pi i\eta}\tilde v(\tau)$, so
$\tilde v(\tau)=0$. 

Let us check that $\tilde v(2\eta)$ vanishes. By using the functional equation
\eqref{e-feO} for $\Omega_{2\eta}$, we obtain
\begin{eqnarray*}
V(2\eta,\mu,\tau,\sigma)&=&c\,
e^{-{\pi i\mu}}
\int \Omega_{2\eta}(t,\tau,\sigma)
\frac{\theta(t+2\eta,\tau)\theta(\mu+t,\sigma)}
       {\theta(t-2\eta,\tau)\theta(4\eta,\tau)
         \theta(t-2\eta,\sigma)\theta(\mu+2\eta,\sigma)}\,dt
\\
&=&c\,
e^{-{\pi i\mu}-4\pi i\eta}
\int \Omega_{2\eta}(t+\sigma,\tau,\sigma)
\frac{\theta(\mu+t,\sigma)}
       {\theta(4\eta,\tau)
         \theta(t-2\eta,\sigma)\theta(\mu+2\eta,\sigma)}\,dt
\\
&=&c\,
e^{-{\pi i\mu}-4\pi i\eta}
\int \Omega_{2\eta}(t,\tau,\sigma)
\frac{\theta(\mu+t-\sigma,\sigma)}
       {\theta(4\eta,\tau)
         \theta(t-2\eta-\sigma,\sigma)\theta(\mu+2\eta,\sigma)}\,dt
\\
&=&c\,
e^{{\pi i\mu}}
\int \Omega_{2\eta}(t,\tau,\sigma)
\frac{\theta(\mu+t,\sigma)}
       {\theta(4\eta,\tau)
         \theta(t-2\eta,\sigma)\theta(\mu+2\eta,\sigma)}\,dt
\\
&=&
\frac{\theta'(0,\tau)}{\theta(4\eta,\tau)}
\,\mathrm{res}_{\lambda=-2\eta}V(\lambda,\mu,\tau,\sigma).
\end{eqnarray*}
In this calculation the change of variable $t\mapsto t-\sigma$
was used. For this our choice of $t$-integration contour is
essential, since it implies that one does not encounter poles when
one deformes it back to the original position.
For general $m$ this identity is part III of Theorem 26 in \cite{FV4}.
Thus 
\[
\tilde v(2\eta)=\frac{\theta'(0,\tau)}{\theta(4\eta,\tau)}
\,\mathrm{res}_{\lambda=-2\eta}\tilde v(\lambda).
\]
 But it follows from
(ii) that
\[
\tilde v(2\eta)=-\frac{\theta'(0,\tau)}{\theta(4\eta,\tau)}
\,\mathrm{res}_{\lambda=-2\eta}\tilde v(\lambda),
\]
so $\tilde v(2\eta)=0$. The same argument may be applied to $2\eta+r+s\tau$
with $r,s\in\{0,\pm 1\}$ (or even for general $r,s$). We have
\[
V(2\eta+r+s\tau,\tau,\sigma)=\frac{\theta'(0,\tau)}{\theta(4\eta,\tau)}
\,e^{2\pi is\sigma}
\,\mathrm{res}_{\lambda=-2\eta+r+s\tau}
V(\lambda,\mu,\tau,\sigma).
\]
This implies that 
\[
\tilde v(2\eta+r+s\tau)=e^{-4\pi i\eta\kappa s}
\frac{\theta'(0,\tau)}{\theta(4\eta,\tau)}\,
\mathrm{res}_{\lambda=-2\eta+r+s\tau}\tilde v(\lambda).
\]
On the other hand, using (ii) and (i) we obtain the same equation but
with the opposite sign, so that both sides vanish. Thus $\tilde v(2\eta+r+s\tau)=0$
and $\tilde v$ is regular at the potential singularities $\lambda=-2\eta+r+s\tau$,
$r,s\in\Z$. $\square$

\medskip

A more direct reformulation of this theorem is the following.

\begin{corollary} 
 Let $m,\kappa\in\Z_{\geq0}$,
$\kappa\geq 2m+2$ and suppose that $\mathrm{Im}\,\eta<0$,
 $\mathrm{Im}\,\tau>0$. Let, for $t\in\C^m$,
\[
\omega_{m}(t,\lambda,\tau)=
\prod_{1\leq i<j\leq m}
\frac{\theta(t_i-t_j,\tau)}
{\theta(t_i-t_j+2\eta,\tau)}
\prod_{j=1}^m
\frac{\theta(\lambda\!+\!t_j,\tau)}
{\theta(t_j-2\eta m,\tau)}\,.
\]
Introduce the integral kernel
\[
M(\lambda,\mu,\tau,p)=\frac{e^{-\frac{\pi i}{4\eta}(\lambda+\mu)^2}
u_0(\lambda,\mu,\tau,p)
\theta(\mu,p)}
{\prod_{j=-m}^m\theta(\lambda-2\eta j,\tau)},
\]
where 
\begin{eqnarray*}
u_0(\lambda,\mu,\tau,p)&=&
\int
\prod_{i=1}^m\Omega_{2\eta m}(t_i,\tau,p)
\prod_{1\leq i<j\leq m}\Omega_{-2\eta}(t_i-t_j,\tau,p)
\\
&&\times
\omega_m(t,\lambda,\tau)\,\omega_m(t,\mu,p)
\,dt_1\cdots dt_m.
\end{eqnarray*}
The integration is over a torus as in \ref{suse-12}.

Then the integral operator
\[
M(\tau)\phi(\lambda)=\int_{2\eta\R}M(\lambda,\mu,\tau,\tau-2\eta\kappa)
\phi(-\mu)\,d\mu
\]
maps $\Theta_{\kappa-2m-2}(\tau-2\eta\kappa)^W$ to
 $\Theta_{\kappa-2m-2}(\tau)^W$.
\end{corollary}

\medskip

\subsection{Remark} A section $v$ of $E_{\kappa,2m,\eta}$ 
is called {\em projectively horizontal} if it obeys the qKZB equation
 $T_{\kappa,m}(\tau)v(\tau-2\eta\kappa)=C(\tau)\,v(\tau)$ up to a
scalar factor $C(\tau)$. For $m=0$ projectively horizontal sections are
given by odd theta functions as in the differential case, see \eqref{e-s1}.
In a sequel \cite{FV3} to this paper, we show that 
for $m=1$ (and conjecturally for higher $m$ as well),
projectively horizontal sections are again given by elliptic hypergeometric
integrals.

\subsection{Remark} 
The compatibility of the difference operator $T_{\kappa,m}(\tau)$ with
the $\mathrm{SL}(2,\Z)$ action can be better understood in terms of
a discrete connection on a space with an $\mathrm{SL}(3,\Z)$-action. This
will be discussed in \cite{FV3}.

\end{document}